\documentclass{ifacconf}
\usepackage{natbib}

\usepackage[utf8]{inputenc}
\usepackage[T1]{fontenc}

\usepackage{amsmath,amsfonts,amssymb} %,amsthm}
\usepackage{mathtools}
\usepackage{mathabx}
\usepackage{mathrsfs}
\usepackage{dsfont}
\usepackage{scalerel}
\usepackage{listings}
\usepackage{graphicx}
\usepackage{enumerate}
\usepackage{url}
\usepackage{csquotes}
\usepackage{xcolor}
\usepackage{booktabs}
\usepackage{multicol}
\usepackage{multirow}
\usepackage{makecell}
\usepackage{soul}

\usepackage{graphicx}
\DeclareGraphicsExtensions{.eps}

\usepackage{calc}
\usepackage{tikz}
\usetikzlibrary{shapes, decorations.markings, calc}
\usetikzlibrary{positioning,hobby,arrows,calc}

\definecolor{blue1}{RGB}{128, 191, 255}
\definecolor{blue2}{RGB}{51   153  255}
\definecolor{blue3}{RGB}{50,60,200}
\definecolor{blue4}{RGB}{30,120,200}
\definecolor{redd1}{RGB}{255, 153, 128}
\definecolor{redd2}{RGB}{255  71   26}
\definecolor{redd3}{RGB}{200,50,50}
\definecolor{gren1}{RGB}{153, 230, 153}
\definecolor{gren2}{RGB}{51   204  51}
\definecolor{gren3}{RGB}{50,140,50}
\definecolor{yelw1}{RGB}{255, 179, 26}
\definecolor{yelw2}{RGB}{255  204  0}
\definecolor{yelw3}{RGB}{241,188,49}
\definecolor{prpl1}{RGB}{204, 153, 255}
\definecolor{prpl2}{RGB}{153  51   255}
\definecolor{mage1}{RGB}{255, 51, 51}
\definecolor{gray1}{RGB}{150,150,150}
\definecolor{brwn1}{RGB}{120,70,20}

\def \S{\mathtt{S}}
\def \I{\mathtt{I}}
\def \D{\mathtt{D}}
\def \H{\mathtt{H}}
\def \E{\mathtt{E}}
\def \R{\mathtt{R}}
\def \d{\text{d}}

\def \ba{\begin{array}}
\def \ea{\end{array}}
\def \bea{\begin{eqnarray}}
\def \eea{\end{eqnarray}}
\def \be{\begin{equation}}
\def \ee{\end{equation}}

\def \colsep{\arraycolsep}

\def \bb{\mathbb}

\def \mc{\mathcal}

\def \diag{\mathtt{diag}}

\def \sym{\mathtt{sym}}

\def \T{\top}%{\intercal}%{\mathtt{T}}

\DeclareMathOperator*{\esssup}{ess\,sup}

\usepackage{textcomp}
\def\BibTeX{{\rm B\Kern-.05em{\sc i\Kern-.025em b}\Kern-.08em
    T\Kern-.1667em\lower.7ex\hbox{E}\Kern-.125emX}}

\begin{document}
\begin{frontmatter}

\title{Feedback Design for Devising Optimal Epidemic Control Policies}

\thanks[footnoteinfo]{This work is supported the European Union's Horizon Research and Innovation Programme under Marie Sk\l{}odowska-Curie grant agreement No. 101062523 and the US National Science Foundation, grant NSF-ECCS \#2032258.}

\author[Umar,Kalle]{Muhammad Umar B. Niazi} 
\author[Phil]{Philip E. Par\'{e}} 
\author[Kalle]{Karl H. Johansson}

\address[Umar]{Laboratory for Information and Decision Systems, Massachusetts Institute of Technology, 77 Massachusetts Avenue, Cambridge, MA 02139, USA (Email: \texttt{niazi@mit.edu})}
\address[Kalle]{Division of Decision and Control Systems, Digital Futures, School of Electrical Engineering and Computer Science, KTH Royal Institute of Technology, Stockholm, Sweden (Email: \texttt{kallej@kth.se})}
\address[Phil]{Elmore Family School of Electrical and Computer Engineering, Purdue University, IN, USA (Email: \texttt{philpare@purdue.edu})}

\begin{abstract}
This paper proposes a feedback design that effectively copes with uncertainties for reliable epidemic monitoring and control. There are several optimization-based methods to estimate the parameters of an epidemic model by utilizing past reported data. However, due to the possibility of noise in the data, the estimated parameters may not be accurate, thereby exacerbating the model uncertainty. To address this issue, we provide an observer design that enables robust state estimation of epidemic processes, even in the presence of uncertain models and noisy measurements. Using the estimated model and state, we then devise optimal control policies by minimizing a predicted cost functional. To demonstrate the effectiveness of our approach, we implement it on a modified SIR epidemic model. The results show that our proposed method is efficient in mitigating the uncertainties that may arise in epidemic monitoring and control.
\end{abstract}


\end{frontmatter}

\section{Introduction}

Since the onset of the SARS-CoV-2 outbreak, there has been a surge of interest in epidemic processes from many fields including the controls community. These works typically consider analysis, parameter identification, state estimation, forecasting, and/or control of a particular compartmental model that may or may not be networked, e.g., \citep{hota2021}. Further, there is also a rich body of literature from the controls field prior to the COVID-19 Pandemic~\citep{wan2007,nowzari2016,mei2017,pare2018}. In this work, we present a unified framework for parameter estimation, state estimation, and optimal control on a generic class of nonlinear models that includes most of the deterministic epidemiological spreading models in the literature. 

The parameter and state estimation problems are questions of identifiability and observability, respectively. Conventionally, differential geometric techniques were employed for obtaining sufficient conditions that verify these notions for nonlinear systems \citep{grewal1976, hermann1977}. 
% An open-source software STRIKE-GOLDD by \cite{villaverde2016} is based on these techniques. 
On the other hand, to obtain necessary and sufficient conditions, \cite{diop1991} introduced differential algebraic methods for identifiability and observability. The concepts were further developed and applied to biological models by \cite{audoly2001} and \cite{saccomani2003}. 
% Software packages DAISY \citep{saccomani2019} and GenSSI \citep{ligon2018} are based on testing these differential algebraic conditions. 

Once identifiability has been verified, the parameters must be estimated, for which several methods exist in the literature. Most common among them are the gradient-based and Newton-type methods like Levenberg–Marquardt and trust region reflective algorithms; see \citep{bard1974} and \citep{ljung1999}. %, and \cite{chowell2017}. 
\cite{barz2015} present some useful regularization techniques to ensure that the parameter estimation problem is well-posed. For some basic epidemic models, like SIR, explicit expressions for parameters are derived by \cite{hadeler2011} and \cite{magal2018}. These ideas have also been explored for networked epidemics~\citep{pare2020}. However, if the state variables in those expressions cannot be directly measured, these techniques cannot be employed. 

On the other hand, constrained optimal control is also a rich area of research with two techniques typically employed in practice. The first one is to use Pontryagin's minimum principle for computing the optimal control trajectory \cite[Chapter 5]{kirk2004}. However, in general, this principle is only a necessary condition of optimality. It is sufficient only in the case when the Hamiltonian functional is convex in the state variable \cite[Chapter 2]{sethi2019}. A practically superior method to solve these types of problems is to convert them to a constrained nonlinear optimization problem \citep{betts2010} and then use a numerical solver. Optimal control has been employed for epidemic mitigation in several works, e.g., \citep{kohler2020, acemoglu2021}.

Unlike parameter estimation and optimal control, the literature on state estimation of epidemic processes is lacking. Designing robust observers to accurately estimate the current state of the epidemic model is the missing component in the feedback optimal control. However, due to being nonlinear, observer design for epidemic processes is quite challenging. Extended Kalman filtering techniques \citep{rajaei2021,gomez2021,azimi2022} are based on linearization and can only provide local guarantees. Therefore, one must know the initial state quite accurately to obtain a good state estimate when using these techniques. On the other hand, observer design techniques for general nonlinear systems with global guarantees turn out to be very conservative for epidemic models \citep{niazi2022cdc}. 

Our main contributions include the extension of the observer proposed by \cite{niazi2022cdc} and provide robust guarantees under model and data uncertainties. Moreover, devising optimal epidemic policies using state estimates from the observer is another novelty considered in this paper. Our results demonstrate that incorporating a robust observer in the feedback loop yields more reliable epidemic control policies.

The rest of the paper is organized as follows. Section~\ref{sec:prob} presents the generic class of nonlinear epidemic models and formulates the problem of interest. 
Section~\ref{sec:outline} outlines our proposed framework and Section~\ref{sec:back} provides a necessary background.
Section~\ref{sec:estimation} presents our proposed algorithm for robust state estimation and Section~\ref{sec:sim} demonstrates our method on a modified SIR epidemic model.

\textit{Notations.} The Euclidean norm of $x\in\bb{R}^n$ is denoted as $\|x\|\doteq \sqrt{x^\T x}$.
For a function $w\in L^\infty(\bb{R};\bb{R}^n)$, the essential supremum norm $\|w\|_{\infty}\doteq\esssup_{t\in\bb{R}} \|w(t)\|$.
By $w_{[t_0,t_1]}$, we denote the restriction of $w$ to $[t_0,t_1]$ for some $t_1>t_0$. 
The maximum singular value of $M\in\bb{R}^{n\times m}$ is denoted as $\sigma_{\max}(M)$. An identity matrix of size $n\times n$ is $I_n$. For $M\in\bb{R}^{n\times n}$, $\sym(M)\doteq M+M^\T$, and $M\geq 0$ ($M>0$) means that $M$ is positive semi-definite (resp., definite).

\section{Problem Definition}\label{sec:prob}

We consider a class of deterministic epidemic models that can be described by
\begin{subequations}
\label{eq:system}
    \begin{align}
        \label{eq:system-state}
        \dot{x}(t) &= Ax(t) + Gf(Hx(t),u(t)) \\
        \label{eq:system-output}
        y(t) &= Cx(t)
    \end{align}
\end{subequations}
where $x(t)\in\mc{X}\subset\bb{R}^{n_x}$ is the state, $u(t)\in\mc{U}\subset\bb{R}^{n_u}$ is the control input, and $y(t)\in\bb{R}^{n_y}$ is the measured output. 
Each element $x_i(t)$ of the state vector corresponds to a different epidemic variable or compartment, and each element $y_i(t)$ of the output vector corresponds to a certain measurable epidemic variable.
% with the zero-order hold. The error incurred by the zero-order hold can also be incorporated in the the measurement noise $v(t)$.
On the other hand, each element $u_i(t)$ of the input vector corresponds to a certain pharmaceutical or non-pharmaceutical interventions like improving medical facilities, testing and isolation, quarantining, vaccination, lockdown, social distancing, and travel restrictions, which can be enforced by a public authority.

We have $A\in\bb{R}^{n_x\times n_x},G\in\bb{R}^{n_x\times n_g},C\in\bb{R}^{n_y\times n_x}$ with
\begin{align*}
    A &\doteq A(\theta), \quad 
    G \doteq G(\theta), \quad
    C \doteq C(\theta)
\end{align*}
where $\theta$ is the vector of epidemic parameters.
The matrix $H\in\{0,1\}^{n_H\times n_x}$ is known and specifies the state variables involved in the nonlinear function $f:\mc{X}\times\mc{U}\rightarrow\bb{R}^{n_f}$, where $f$ is smooth and thus Lipschitz continuous on a compact domain $\mc{X}\times\mc{U}$. That is, for every $x,\hat{x}\in\mc{X}$ and $u\in\mc{U}$, there exists $\ell\geq 0$ such that
\be \label{eq:lipschitz_property}
\|f(Hx,u)-f(H\hat{x},u)\| \leq \ell \|Hx-H\hat{x}\|
\ee
where
\be \label{eq:lipschitz_constant}
\ell = \sup_{(x,u)\in\mc{X}\times\mc{U}} \sigma_{\max} \left(\frac{\partial f}{\partial x}(Hx,u) \right).
\ee
Note that $f$ depends only on the state $x(t)$ and the input $u(t)$, and not on the parameters $\theta$.
% Note that \eqref{eq:lipschitz_constant} is a constrained nonlinear optimization problem that can be solved by, for instance, the \texttt{fmincon} solver in MATLAB.

\begin{rem}
    The class of nonlinear systems \eqref{eq:system} captures a variety of epidemic models in the literature. For instance, all the basic SIS, SIR, SEIR models \citep{hethcote1989, hethcote1994, hethcote2000} and their variants \citep{arino2003,giordano2020,giordano2021,niazi2021arc} can be written in the form of \eqref{eq:system-state}. The networked epidemic models \citep{nowzari2016,mei2017,pare2020} can also written compactly as \eqref{eq:system-state}. \quad $\triangle$
\end{rem}

The reported data on a time interval $[t_0,t_1]$ is given by
\begin{subequations}
\label{eq:input-output-data}
    \begin{align}
        \bar{u}(t) &= u(t) + \delta_u(t) \label{eq:input-output-data-input} \\
        \bar{y}(t) &= y(t) + \delta_y(t) \label{eq:input-output-data-output}
    \end{align}
\end{subequations}
where $t_0\geq 0$ is the time of epidemic onset, $t_1>t_0$ is the current time, and $\delta_u(t)$ and $\delta_y(t)$ represent the uncertainties in the input-output data. The uncertainties $\delta_u(t)$ and $\delta_y(t)$ are unknown and account for clerical errors and delays in recording and reporting the data.

\textit{Problem statement.} Given the input-output data $(\bar{u},\bar{y})$ for the past time interval $[t_0,t_1]$, we first aim to estimate the parameters $\theta$ and the current state $x(t_1)$ of \eqref{eq:system}. Then, based on the estimated model, we devise an optimal control policy $u(t)$ for a future time interval $[t_1,t_2]$, $t_2>t_1$, by minimizing a given cost functional 
\be \label{eq:cost_fun}
    J(x_{[t_1,t_2]},u_{[t_1,t_2]}) \doteq  \int_{t_1}^{t_2} q(x,u,t) \d t
\ee
subject to a set of specified constraints
\be \label{eq:constraints}
\ba{rclrcl}
    r_i(x_{[t_1,t_2]},u_{[t_1,t_2]}) &=& 0, & \quad i &=& 1,2,\dots,k \\
    s_j(x_{[t_1,t_2]},u_{[t_1,t_2]}) &\leq& 0, & \quad j &=& 1,2,\dots,l
\ea
\ee
where
$x_{[t_1,t_2]}\doteq x(t, u; \hat{x}_{t_1}, \hat{\theta} )$
is the predicted state trajectory obtained by integrating \eqref{eq:system-state} for $t\in[t_1,t_2]$ using the values of estimated parameters $\theta=\hat{\theta}$ and choosing the initial condition as the estimated state $x(t_1)=\hat{x}_{t_1}$.

% \begin{rem}
%     The problem state above can also be formulated as a receding horizon optimal control problem. However, for simplicity, we consider an optimal control  
% \end{rem}

% The optimal control policy $u$ minimizes the predicted cost
% \[
%     J(x_{[t_1,t_2]},u_{[t_1,t_2]}) \doteq  \int_{t_1}^{t_2} q(x,u,t) \d t
% \]
% subject to a set of constraints, where $q:\mc{X}\times\mc{U}\times\bb{R}_{\geq 0}\rightarrow\bb{R}_{\geq 0}$ is a smooth and convex functional, and 
% \[
% x_{[t_1,t_2]}\doteq x_{[t_1,t_2]}(t, u; \hat{x}_{t_1}, \hat{\theta} )
% \]
% is the predicted state trajectory obtained by integrating \eqref{eq:system-state} for $t\in[t_1,t_2]$ using the values of estimated parameters $\theta=\hat{\theta}$ and initial condition $x(t_1)=\hat{x}_{t_1}$.

\section{Outline of the Proposed Method}\label{sec:outline}

Given the past input-output data \eqref{eq:input-output-data} and the cost functional \eqref{eq:cost_fun} with constraints \eqref{eq:constraints}, the proposed feedback design for \eqref{eq:system} has three main constituents.

\begin{enumerate}
\itemsep0.5em
    \item \textbf{Parameter estimation.} Given the past input-output data $(\bar{u}(t),\bar{y}(t))$, for $t\in[t_0,t_1]$, we estimate the model parameters $\theta$ by solving 
    \be \label{prob:par_id}
        \hat{\theta} = \arg\min_{\theta\in\Theta} \int_{t_0}^{t_1} \|\bar{y}(t) - y(t,\bar{u}; \theta) \| \d t
    \ee
    where $y(t,\bar{u};\theta)$ is the output of \eqref{eq:system} at time $t$ given the parameters $\theta$.
    
    \item \textbf{State estimation.} Given the past input-output data $(\bar{u}(t),\bar{y}(t))$, for $t\in[t_0,t_1]$, we design a state observer that estimates the current state $x(t_1)$ by solving
    \begin{equation}
    \label{prob:state_est}
        \hat{x}(t_1) = \arg\min_{x_{t_1}\in\mc{X}} \int_{t_0}^{t_1} \|\bar{y}(t) - y(t,\bar{u};x_{t_1},\hat{\theta})\| \d t
    \end{equation}
    where $y(t,\bar{u};x_{t_1},\hat{\theta})$ is the output of \eqref{eq:system} at time $t$ when the state trajectory goes through $x_{t_1}$ at time $t_1$ and the parameters $\theta=\hat{\theta}$.
        
    \item \textbf{Optimal control:} Given the estimated model \eqref{eq:system} with parameters $\hat{\theta}$ and the state $\hat{x}(t_1)$, we obtain optimal control policies by solving
    \begin{subequations}
    \label{prob:opt_con}
        \begin{align}
            u^*(t) = & \arg\min_{u\in\mc{U}} J(x,u) \\
            & \text{subject to \eqref{eq:constraints} and}~ \nonumber \\ 
            & \left\{\colsep=0pt\ba{ccl} 
            \dot{x} &=& A(\hat{\theta})x + G(\hat{\theta}) f(Hx,u) \\ [0.5em]
            t &\in & [t_1,t_2], \; x(t_1) = \hat{x}(t_1)
            \ea\right. \nonumber
        \end{align}
    \end{subequations}
    where the cost functional $J(x,u)$ is defined in \eqref{eq:cost_fun}.
\end{enumerate}

\begin{figure}[!]
    \begin{center}
        \begin{tikzpicture}[scale=0.8]
            \node[rectangle, draw, text width=2.5cm, align=center, fill=green!50!blue!20] (OC) at (0,4) {Optimal Control};
            \node[rectangle, draw, text width=2.5cm, align=center, fill=yellow!50!red!30] (EP) at (4,4) {Epidemic Process};
            
            \node[rectangle, draw, text width=2.5cm, align=center, fill=red!40!green!20] (SE) at (2,-0.25) {State\\ Estimation};
            \node[rectangle, draw, text width=2.5cm, align=center, fill=blue!50!red!20] (PE) at (2,1.5) {Parameter Estimation};
            
            \node[circle, draw, fill=red!50!black!5] (sum1) at (6,2.75) {\scriptsize $+$};
            \node[circle, draw, fill=red!50!black!5] (sum2) at (1.9,2.75) {\scriptsize $+$};
            
            \node[anchor=north] at (6,1.75) {\small $\delta_y$};
            \node[anchor=east] at (0.9,2.75) {\small $\delta_u$};

            \draw[fill=black!60,black!60] (4.2,2.25) rectangle (4.8,2.35);
            \draw[fill=black!60,black!60] (-1.8,2.45) rectangle (-2.4,2.55);
            
            \draw[-latex,thick,black!60] (1.9,4) -- (sum2);
            \draw[-latex,thick,black!60] (sum1) -- node[pos=0.9,above] {\small $\bar{y}$} (4.7,2.75) -- (4.7,2.35);
            \draw[-latex,thick,black!60] (sum2) -- node[pos=0.9,above] {\small $\bar{u}$} (4.3,2.75) -- (4.3,2.35);
            \draw[-latex,thick,black!60] (4.5,2.25) -- (4.5,-0.25) -- (SE);
            \draw[-latex,thick,black!60] (4.5,1.5) -- (PE);

            \draw[-latex,thick,black!60] (6,1.75) -- (sum1);
            \draw[-latex,thick,black!60] (0.9,2.75) -- (sum2);
            
            \draw[-latex,thick,black!60] (OC) -- node[pos=0.5,above] {\small $u$} (EP);
            \draw[-latex,thick,black!60] (EP) -- node[pos=1.25,above] {\small $y$} (6,4) -- (sum1);

            \draw[-latex,thick,black!60] (PE) -- node[pos=0.2,above] {\small $\hat{\theta}$} (-1.9,1.5) -- (-1.9,2.45);
            \draw[-latex,thick,black!60] (SE) -- node[pos=0.175,above] {\small $\hat{x}$} (-2.3,-0.25) -- (-2.3,2.45);
            \draw[-latex,thick,black!60] (-2.1,2.55) -- (-2.1,4) -- (OC);

            \draw[-latex,thick,black!60] (PE) -- node[pos=0.5,right] {\small $\hat{\theta}$} (SE);
        \end{tikzpicture}
    \end{center}
    \caption{Block scheme of optimal feedback control.}
    \label{fig:algo_block}
\end{figure}
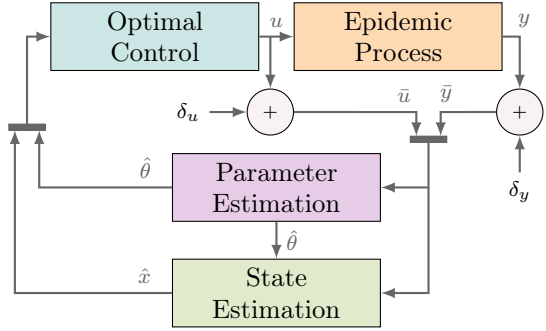

\section{Background Material}\label{sec:back}

In this section, we introduce notions of identifiability and observability, and describe techniques to solve the parameter estimation and optimal control problems.

\subsection{Verifying identifiability and observability}

Identifiability is necessary for estimating model parameters from the input-output data $(\bar{u},\bar{y})$. If the model has non-identifiable parameters, then it is not possible to estimate them uniquely \citep{saccomani2003}.

\begin{defn}
    System \eqref{eq:system} is \textit{locally identifiable} if, for almost any parameter vector $\theta\in\Theta$, there exists a neighborhood $\mc{N}(\theta)\subseteq\Theta$ such that, for every $\hat{\theta}\in\mc{N}(\theta)$,
    \[
        y(t;\theta)=y(t;\hat{\theta}), \forall t\geq 0  \Leftrightarrow \theta=\hat{\theta}.
    \]
\end{defn}

On the other hand, the notion of observability guarantees whether or not input-output data contains sufficient information about the system's state \citep{bernard2022}.

\begin{defn}
    System \eqref{eq:system} is \textit{locally observable} if, for any $\tau\geq 0$ and almost any $x(\tau)\in\mc{X}$, there exists a neighborhood $\mc{N}(x(\tau))\in\mc{X}$ such that, for every $\hat{x}(\tau)\in\mc{N}(x(\tau))$,
    \[
        y(t,u;x(\tau)) = y(t,u;\hat{x}(\tau)), \forall t\in[\tau,\infty) \Leftrightarrow x(\tau) = \hat{x}(\tau).
    \]
\end{defn}

Identifiability and observability of \eqref{eq:system} are required for the well-posedness of the parameter and state estimation problems, respectively.
To verify these notions for \eqref{eq:system}, we employ the GenSSI algorithm proposed by \cite{ligon2018}, which tests the injectivity of the observation map obtained by taking Lie derivatives of the output.

\subsection{Solving the parameter estimation problem}
\label{subsec_parest}

Once identifiability of the system has been verified, the parameters must be estimated. %In order t
To carry out this process, optimization algorithms like Levenberg-Marquardt \citep{more1978} and trust region \citep{byrd1987} can be used. 
%Such algorithms are already integrated in the System Identification Toolbox of MATLAB, which is used in this paper for parameter estimation.

\subsection{Solving the optimal control problem}

In order to solve \eqref{prob:opt_con}, we convert it to a constrained nonlinear optimization problem \citep{betts2010} and use interior point \citep{byrd1999} or trust region reflective \citep{coleman1994} methods. 
%We propose using optimization solvers that implement these algorithms, like \texttt{fmincon} in MATLAB and CasADi \citep{andersson2019}. 

\section{Algorithm for Robust State Estimation}\label{sec:estimation}

The main challenge in the proposed scheme (Figure~\ref{fig:algo_block}) is the state estimation problem. Existing observer design techniques for the state estimation of epidemic processes are quite conservative and often turn out to be infeasible \citep{niazi2022cdc}. Moreover, observers for epidemic models are often designed for specific compartmental models and cannot be adapted to other models. Here, we extend an observer for general epidemic models proposed by \cite{niazi2022cdc} to include model uncertainties and measurement noise. 
%We also provide a sufficient condition for the robustness of the observer. % to parameter and data uncertainties.

After estimating the parameters from the data \eqref{eq:input-output-data}, \eqref{eq:system} can be written as an uncertain nonlinear system
\begin{subequations}
\label{eq:system_uncertain}
    \begin{align}
        \dot{x}(t) &= \hat{A}x(t) + \hat{G}f(Hx(t),\bar{u}(t)) + w(t) \\
        \bar{y}(t) &= \hat{C}x(t) + v(t)
    \end{align}
\end{subequations}
where $w(t)\in\bb{R}^{n_x}$ is the model uncertainty, $v(t)\in\bb{R}^{n_y}$ is the measurement noise,
% . Also, in \eqref{eq:system_uncertain}, we have
and 
\[
\hat{A} \doteq A(\hat{\theta}), \quad \hat{G} \doteq G(\hat{\theta}), \quad \hat{C} \doteq C(\hat{\theta}).
\]
Notice that the model uncertainty and measurement noise result from the uncertainties in the input-output data and the parameter estimation error $\theta-\hat{\theta}$. 

Consider the observer proposed by \cite{niazi2022cdc}:
\begin{subequations}
\label{eq:observer}
    \begin{align}
        \dot{z}(t) &= Mz(t) + (ML+J)\bar{y}(t) + N\hat{G}f(q(t),\bar{u}(t)) \\
        \hat{x}(t) &= z(t) + L\bar{y}(t) \\
        \hat{y}(t) &= \hat{C}\hat{x}(t)
    \end{align}
\end{subequations}
where $q(t) \doteq H\hat{x}(t) + K(\bar{y}(t)-\hat{y}(t))$, $J,L\in\bb{R}^{n_x\times n_y}$ and $K\in\bb{R}^{n_H\times n_y}$ are matrices to be designed, and
\begin{align*}
    M = \hat{A}-L\hat{C}\hat{A} -J\hat{C}, \quad
    N = I_{n_x} - L\hat{C}.
\end{align*}
Here, $z(t)\in\bb{R}^{n_x}$ is the observer's state, and $\hat{x}(t)\in\bb{R}^{n_x}$ and $\hat{y}(t)\in\bb{R}^{n_y}$ are the state and output estimate of \eqref{eq:system_uncertain}.

Consider the semidefinite programming (SDP) problem:
\begin{subequations}
\label{eq:sdp}
    \begin{align}
        &
        \text{minimize}~ \mu ~\text{subject to}\\
        &
        \left[\ba{cc}
        \sym(P\hat{A} - R\hat{C}\hat{A} - S\hat{C}) + Q & (P-R\hat{C})\hat{G} \\
        G^\T(P-R\hat{C})^\T & -I_{n_f}
        \ea\right] < 0  \label{eq:lmi1} \\
        &
        \left[\ba{cc}
        -Q & (H-K\hat{C})^\T \\
        H-K\hat{C} & -\frac{1}{\ell^2} I_{n_H}
        \ea\right] \leq 0 \label{eq:lmi2} \\
        &
        \left[\ba{cc}
        -\mu I_{n_x} & R \\ R^\T & - I_{n_y}
        \ea\right] \leq 0 \label{eq:lmi3} \\
        & 
        P=P^\T>0 ~\text{and}~ Q=Q^\T>0 
    \end{align}
\end{subequations}
where $\ell$ is the Lipschitz constant obtained from \eqref{eq:lipschitz_constant}, and
\[J=P^{-1} S, \quad L=P^{-1} R.\]
Now, we show that feasibility of \eqref{eq:sdp} is sufficient for the existence of a robust observer \eqref{eq:observer}.

\begin{thm}
\label{thm:observer}
    If the SDP problem \eqref{eq:sdp} is feasible, then there exist a $\mc{KL}$ function $\beta$ and $\mc{K}_\infty$ functions $\alpha_1,\alpha_2$ such that the estimation error satisfies
    \begin{multline*}
        \|x(t_1)-\hat{x}(t_1)\| \leq \beta(\|x(t_0)-\hat{x}(t_0)\|,t_1) \\ + \alpha_1(\|w_{[t_0,t_1]}\|_{\infty}) + \alpha_2(\|v_{[t_0,t_1]}\|_{\infty}).
    \end{multline*}
\end{thm}
\begin{pf}[Sketch]
    Define the estimation error 
    \be \label{eq:estimation_error}
    e \doteq x-\hat{x}= x - (z+L\bar{y}) = \eta-Lv
    \ee
    where $x,\bar{y}$ are given in \eqref{eq:system_uncertain}, $z,\hat{x}$ are given in \eqref{eq:observer}, and
    $
    \eta \doteq Nx-z.
    $
    Then, by taking the time derivative of $\eta$, we obtain an error system
    \begin{align*}
        \dot{\eta} &= M\eta + NG \tilde{f} + Nw - (ML+J) v, \quad 
        e = \eta - Lv
    \end{align*}
    where $\tilde{f}\doteq f(Hx,\bar{u})-f(H\hat{x}+K(\bar{y}-\hat{C}\hat{x}),\bar{u})$ and the dependence on $t$ is omitted for brevity.
    Here, $\eta(t)$ is the state, $e(t)$ is the output, and $w(t),v(t)$ are unknown inputs.
    Finally, subject to the feasibility of \eqref{eq:sdp}, the proof is completed by showing input-to-output stability \citep{sontag2000} of the above system. \qed
\end{pf}

The above theorem states that the state estimation error \eqref{eq:estimation_error} is stable with respect to the data and parameter uncertainties. Particularly, the noise attenuation is directly related to the parameter $\mu$ in \eqref{eq:sdp}. Moreover, in the absence of uncertainties, the estimate $\hat{x}(t_1)$ asymptotically converges to the true $x(t_1)$ given that $t_1$ is sufficiently large. 
As a result of the $\mc{KL}$ function $\beta$, the transient error resulting from the poor choice of $\hat{x}(t_0)$ also converges to zero asymptotically.

%\begin{rem}
%    For the feasibility of \eqref{eq:lmi1}, it is necessary that $(\hat{A},\hat{C})$ be a detectable pair. In case this does not hold, we can always add terms in $\hat{A}$ and subtract them in $G$ by adding extra linear variables in $f(Hx,u)$ so that $(\hat{A},\hat{C})$ becomes detectable. However, this comes at the expense of increasing the Lipschitz constant $\ell$. \quad $\triangle$
%\end{rem}

\section{Application to SIDHER Epidemic Model}\label{sec:sim}

In this paper, we demonstrate the proposed method on an SIDHER epidemic model. After providing the model design, we implement the proposed method step-by-step and provide the simulation results.

\subsection{Model design}

\begin{figure}
    \begin{center}
    \begin{tikzpicture}[scale=0.8]
        \node[rounded corners, draw, text width=0.5cm, minimum height=0.75cm, text centered, fill=red!30!blue!20] (S) at (0,2) {$\S$};
        \node[rounded corners, draw, text width=0.5cm, minimum height=0.75cm, text centered, fill=brown!30!red!30] (I) at (2,2) {$\I$};
        \node[rounded corners, draw, text width=0.5cm, minimum height=0.75cm, text centered, fill=brown!30!red!30] (D) at (4,2) {$\D$};
        \node[rounded corners, draw, text width=0.5cm, minimum height=0.75cm, text centered, fill=brown!30!red!30] (H) at (4,0) {$\H$};
        \node[rounded corners, draw, text width=0.5cm, minimum height=0.75cm, text centered, fill=yellow!30!red!50] (E) at (6,0) {$\E$};
        \node[rounded corners, draw, text width=0.5cm, minimum height=0.75cm, text centered, fill=blue!30!green!30] (R) at (2,0) {$\R$};
        
        \draw[-latex,black!60] (S) to node[pos=0.5, above] {$\beta$} (I);
        \draw[-latex,black!60] (S) to[bend left=10]  node[pos=0.5, above] {$\nu$} (R);
        \draw[-latex,black!60] (I) to node[pos=0.5, above] {$\tau$} (D);
        \draw[-latex,black!60] (I) to node[pos=0.5, right] {$\gamma$} (R);
        \draw[-latex,black!60] (D) to node[pos=0.5, right] {$\phi$} (H);
        \draw[-latex,black!60] (D) to node[pos=0.5, above] {$\rho$} (R);
        \draw[-latex,black!60] (H) to node[pos=0.5, below] {$\sigma$} (R);
        \draw[-latex,black!60] (H) to node[pos=0.5, below] {$\xi$} (E);
        \draw[-latex,black!60] (R) to[bend left=10] node[pos=0.5, below] {$\lambda$} (S);
    \end{tikzpicture}
    \end{center}
    \caption{Block diagram of SIDHER epidemic model.}
    \label{fig:sidher}
\end{figure}
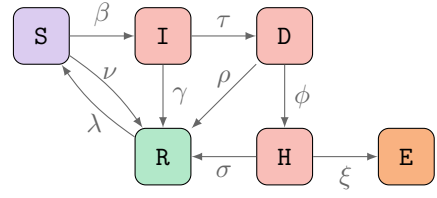

We consider an SIDHER epidemic model (Susceptible, Infected, Detected, Hospitalized, Extinct, and Recovered), which is illustrated in Figure~\ref{fig:sidher} and given by
\begin{subequations}
\label{eq:sidher}
    \begin{align}
        \dot{\S}(t) &= \lambda\R(t) - \beta\S(t)\I(t) (1-u_1(t)) - \nu\S(t) u_4(t) \\
        \dot{\I}(t) &= -\gamma\I(t) + \beta\S(t)\I(t) (1-u_1(t)) - \tau\I(t) u_3(t) \\
        \dot{\D}(t) &= -(\rho+\phi)\D(t) + \tau\I(t) u_3(t) \\
        \dot{\H}(t) &= -\xi\H(t) (1-u_2(t)) + \phi\D(t) - \sigma\H(t) u_2(t) \\
        \dot{\E}(t) &= \xi\H(t) (1-u_2(t)) \\
        \dot{\R}(t) &= -\lambda\R(t) + \gamma\I(t) + \rho\D(t) + \nu\S(t) u_4(t) + \sigma\H(t) u_2(t)
    \end{align}
\end{subequations}
where all the state variables $\S(t),\I(t),\D(t),\H(t),\E(t),\R(t)\in[0,1]$ and control input $u(t)=[\ba{cccc} u_1(t) & u_2(t) & u_3(t) & u_4(t) \ea]^\T$. Note that, for every $t\geq 0$,
\be \label{eq:total_population}
    \S(t) + \I(t) + \D(t) + \H(t) + \E(t) + \R(t) = 1.
\ee

\begin{table}[!]
    \centering
    \caption{Description of control inputs, which are unitless.}
    \label{tab:control}
    \begin{tabular}{c|l}
        \toprule
        \textbf{Control} & \textbf{Description} \\
        \midrule
        $u_1$ & Stringency of NPIs \\
        $u_2$ & Proportion of medical resources dedicated \\
        $u_3$ & Testing capacity per population \\
        $u_4$ & Vaccination capacity per population \\
        \bottomrule
    \end{tabular}
\end{table}
\begin{table}[!]
    \centering
    \caption{Description of parameters.}
    \label{tab:parameters}
    \begin{tabular}{c|l|c}
        \toprule
        \textbf{Parameter} & \textbf{Description} & \textbf{Unit} \\
        \midrule
        $\beta$ & Infection rate & 1/day \\
        $\gamma$ & Recovery rate of undetected cases & 1/day \\
        $\rho$ & Recovery rate of detected cases & 1/day \\
        $\sigma$ & Recovery rate of hospitalized cases & 1/day \\
        $\xi$ & Mortality rate of hospitalized cases & 1/day \\
        $\lambda$ & Rate at which people lose immunity & 1/day \\
        $\phi$ & Hospitalization rate & 1/day \\
        $\tau$ & Testing rate of undetected cases & 1/day \\
        $\nu$ & Vaccination rate of susceptible cases & 1/day \\
        \bottomrule
    \end{tabular}
\end{table}

\textit{Control inputs.}
The control inputs, summarized in Table~\ref{tab:control}, are described below:
\begin{itemize}
    \item $u_1(t)\in[0,1]$, which is the stringency of the non-pharmaceutical interventions (NPIs) like lockdown, social distancing, travel restriction, etc.

    \item $u_2(t)\in[0,1]$, which is a measure of resources invested in improving medical resources and infrastructure, implementing better treatment methods of hospitalized cases, and medical facilities dedicated to the epidemic.

    \item $u_3(t)\in[0,1]$, which is the testing capacity in terms of the fraction of population that can be tested.

    \item $u_4(t)\in[0,1]$, which is the vaccination capacity in terms of the fraction of population that can be vaccinated.
\end{itemize}

\textit{Parameters.}
Summarized in Table~\ref{tab:parameters}, the parameters are described below:
\begin{itemize}
    \item The \textit{infection rate} $\beta\geq 0$ is the product of \textit{contact rate} (average number of contacts each person makes per day) and the \textit{infection probability} (probability that a susceptible person gets infected after coming in contact with an infected person).

    \item The \textit{recovery rates} $\gamma, \rho, \sigma\in[0,1]$ are the inverses of average infection periods of infected, detected, and hospitalized cases, respectively. The \textit{mortality rate} $\xi\in[0,1]$ is the inverse of average number of days after which a typical non-surviving hospitalized case dies. The rates $\lambda,\phi\in[0,1]$ are the inverses of the average number of days after which recovered or vaccinated cases lose their immunity and detected cases are hospitalized, respectively.

    \item The \textit{testing} and \textit{vaccination rates} $\tau,\nu\in[0,1]$ are the fractions of infected individuals that get tested and susceptible individuals that get vaccinated every day on average, respectively.
\end{itemize}

\textit{Measured outputs.}
The model outputs are the following:
\begin{itemize}
    \item $y_1 = \nu\S$: Since we measure the proportion of population vaccinated per day $\nu\S u_4$ and the vaccination capacity $u_4$ is known, we obtain the output $y_1$ by dividing both.    
    \item $y_2 = \tau\I$: Since we measure the proportion of population tested per day $\tau\I u_3$ and the testing capacity $u_3$ is known, we obtain the output $y_2$ by dividing both.
    \item $y_3 = \D$: Active number of detected infected cases.
    \item $y_4 = \rho\D$: Daily number of cases recovering after being detected. This may not be directly measured if people are not tested again after their infectious period. In such a case, this output variable can be inferred from the daily reported data on active detected cases $y_3$, active hospitalized cases $y_6$, and daily number of people recovering from hospitals $y_7$ \citep{niazi2021arc}. 
    \item $y_5 = \phi\D$: Daily number of cases hospitalized after being detected.
    \item $y_6 = \H$: Active number of hospitalized cases.
    \item $y_7 = \sigma\H$: Daily number of cases recovering after hospitalization.
    \item $y_8 = \xi\H$: Daily number of deaths. We assume that non-surviving cases are first detected and then hospitalized as their symptoms worsen.
    \item $y_9 = \E$: Total number of deaths.
    \item $y_{10} = \S + \I + \R$: Since \eqref{eq:total_population} holds and $\D,\H,\E$ are measured, we can obtain the output $y_{10}$ by subtracting $\D+\H+\E$ from $1$.
\end{itemize}

\textit{Model in vector form.}
We can write the model \eqref{eq:sidher} in the form \eqref{eq:system} with $f(Hx,u) = [\ba{ccccc} \S\I & \S\I u_1 & \H u_2 & \I u_3 & \S u_4 \ea]^\T$ and
\begin{align*}
    A &= {\scriptsize \left[\ba{rrcrrr}
    0 & 0 & 0 & 0 & 0 & \lambda \\
    0 & -\gamma & 0 & 0 & 0 & 0 \\
    0 & 0 & -(\rho+\phi) & 0 & 0 & 0 \\
    0 & 0 & \phi & -\xi & 0 & 0 \\
    0 & 0 & 0 & \xi & 0 & 0\\
    0 & \gamma & \rho & 0 & 0 & -\lambda
    \ea\right]}, 
    \quad
    H = {\scriptsize \left[\ba{cccccc}
    1 & 0 & 0 & 0 & 0 & 0 \\
    0 & 1 & 0 & 0 & 0 & 0 \\
    0 & 0 & 0 & 1 & 0 & 0
    \ea\right]},
    \\
    G &= {\scriptsize \left[\ba{rrrrr}
    -\beta & \beta & 0 & 0 & -\nu \\
    \beta & -\beta & 0 & -\tau & 0 \\
    0 & 0 & 0 & \tau & 0 \\
    0 & 0 & -\sigma+\xi & 0 & 0 \\
    0 & 0 & -\xi & 0 & 0 \\
    0 & 0 & \sigma & 0 & \nu
    \ea\right]},
    \quad
    C = {\scriptsize \left[\ba{cccccc}
    \nu & 0 & 0 & 0 & 0 & 0 \\
    0 & \tau & 0 & 0 & 0 & 0 \\
    0 & 0 & 1 & 0 & 0 & 0 \\
    0 & 0 & \rho & 0 & 0 & 0 \\
    0 & 0 & \phi & 0 & 0 & 0 \\
    0 & 0 & 0 & 1 & 0 & 0 \\
    0 & 0 & 0 & \sigma & 0 & 0 \\
    0 & 0 & 0 & \xi & 0 & 0 \\
    0 & 0 & 0 & 0 & 1 & 0 \\
    1 & 1 & 0 & 0 & 0 & 1
    \ea\right]}.
\end{align*}
% From this, notice 
Note that $(A,C)$ is an observable pair.

\begin{figure}[!htp]
    \centering
    \includegraphics[width=0.3\textwidth]{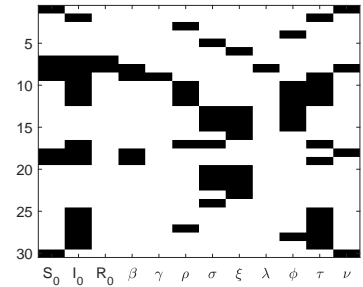}
    \caption{Observability and identifiability tableau.}
    \label{fig:tableau}
\end{figure}

\textit{Identifiability and observability.}
Identifiability and observability of \eqref{eq:sidher} is verified by the GenSSI software \citep{ligon2018} in MATLAB. Fig.~\ref{fig:tableau} illustrates the resulting observability and identifiability tableau, which is the zero-pattern structure of the Jacobian of the observability and identifiability map. The black boxes in the figure represent non-zero terms of the Jacobian, whereas white area represents zero terms. From the figure, it can be seen that the Jacobian has full generic rank, thus implying structural identifiability and observability. 
Therefore, it can be inferred that for almost all values of the initial states and parameters, \eqref{eq:sidher} is at least locally identifiable and observable.

\subsection{Simulation results}

\textit{Data generation.}
We generate noisy synthetic data by simulating the model \eqref{eq:sidher} for $t\in[t_0,t_1]$, where $t_0=0$ and $t_1=30$ days, with ``true'' parameters provided in Table~\ref{tab:true_est_parameters} to illustrate the proposed method. The input-output data $(\bar{u},\bar{y})$ is corrupted with white Gaussian noise that is sampled from $\mc{N}(0,10^{-6})$. The nominal input is chosen as
\[
u(t) = \left[\ba{c} 0.01\lceil \sin(t/2) \rfloor + 0.015 \\ 0.01\lceil \cos(t/2) \rfloor + 0.015 \\ 0.01\lceil \sin(t/3) \rfloor + 0.015 \\ 0.01\lceil \cos(t/3) \rfloor + 0.015 \ea\right].
\]
The true initial state is chosen to be 
\[
x_{t_0}=[\ba{cccccc} 0.999 & 0.0005 & 0.0005 & 0 & 0 & 0 \ea]^\T
\]
which is only used for data generation and is not known by the parameter and state estimation algorithms.

\textit{Parameter estimation.}
The parameters $\rho,\phi,\sigma,\xi$ can be estimated using the least square solution by
\begin{align*}
    \hat{\rho} &= \int_{t_0}^{t_1} \frac{\bar{y}_3(t)^\T \bar{y}_4(t)}{\bar{y}_3(t)^\T \bar{y}_3(t)} \d t, & \hat{\phi} &= \int_{t_0}^{t_1} \frac{\bar{y}_3(t)^\T \bar{y}_5(t)}{\bar{y}_3(t)^\T \bar{y}_3(t)} \d t \\
    \hat{\sigma} &= \int_{t_0}^{t_1} \frac{\bar{y}_6(t)^\T \bar{y}_7(t)}{\bar{y}_6(t)^\T \bar{y}_6(t)} \d t, & \hat{\xi} &= \int_{t_0}^{t_1} \frac{\bar{y}_6(t)^\T \bar{y}_8(t)}{\bar{y}_6(t)^\T \bar{y}_6(t)} \d t.
\end{align*}
The data $\bar{y}$ is discrete and can be interpolated for computing the above integrals.
By fixing these estimated parameters, the remaining parameters are estimated using the \texttt{nlgreyest} function from the System Identification Toolbox in MATLAB employing the Trust-Region-Reflective Algorithm. See Table~\ref{tab:true_est_parameters} for the estimated parameters. Notice that, for parameter estimation, we do not know the true initial state $x_{t_0}$. Instead, we guess the initial state appropriately, where $\S_{t_0}$ is chosen uniformly at random from $[0.95,1]$; $\I_{t_0}$ from $[0,0.05]$; $\D_{t_0},\H_{t_0},\E_{t_0}$ are obtained from $y_3(t_0),y_6(t_0),y_9(t_0)$; and $\R_{t_0}$ is obtained from \eqref{eq:total_population}.

\begin{table}[h]
    \centering
    \caption{True and estimated values of the model parameters.}
    \begin{tabular}{c|c|c}
        \toprule
        \textbf{Parameter} & \textbf{True value} & \textbf{Estimated value} \\
        \midrule
        $\beta$ & 0.3500  &  0.3530 \\
        $\gamma$ & 0.1000  &  0.0981 \\
        $\rho$ & 0.0500  &  0.0501 \\
        $\sigma$ & 0.0400  &  0.0399 \\
        $\xi$ & 0.0200  &  0.0202 \\
        $\lambda$ & 0.0167  &  0.0383 \\
        $\phi$ & 0.1429  &  0.1428 \\
        $\tau$ & 0.3000  &  0.2757 \\
        $\nu$ & 0.0100  &  0.0100 \\
        \bottomrule
    \end{tabular}
    \label{tab:true_est_parameters}
\end{table}

\begin{rem}
    Since the system \eqref{eq:sidher} is locally identifiable, and the parameter estimation problem \eqref{prob:par_id} is non-convex, the solution obtained in Table~\ref{tab:true_est_parameters} might be local. That is, the solution depends on the appropriate initialization of these parameters for the optimization algorithm. To address this issue, one can adopt global optimization techniques to ensure that the true parameters can be efficiently estimated. Moreover, the observer can be extended to simultaneously estimate the state and parameters of the model. However, this is an interesting prospect that will be addressed in the future.  \quad $\triangle$
\end{rem}

\begin{figure}[!thp]
    \centering
    \includegraphics[width=0.475\textwidth]{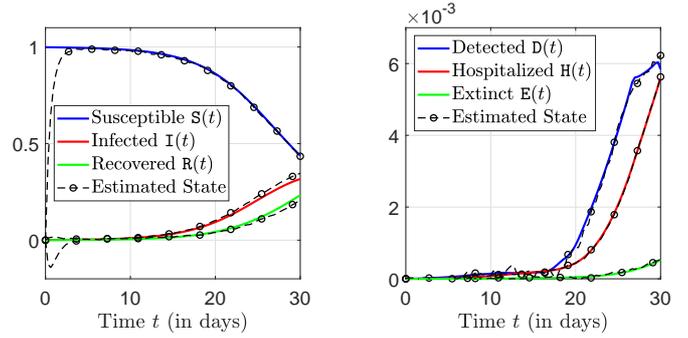}
    \caption{State estimation by the proposed observer.}
    \label{fig:state_est}
\end{figure}

\begin{figure}[!thp]
    \centering
    \includegraphics[width=0.475\textwidth]{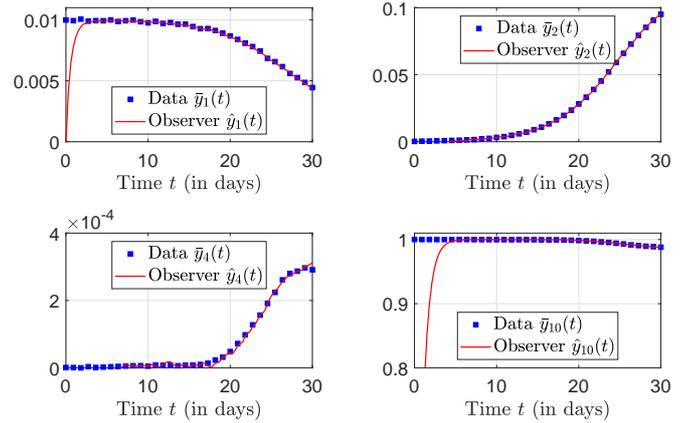}
    \caption{Online output tracking by the observer.}
    \label{fig:pred_y}
\end{figure}

\textit{State estimation.}
After obtaining the estimated parameters, we employ the observer \eqref{eq:observer} for estimating the state of \eqref{eq:sidher}. The observer is designed by solving \eqref{eq:sdp}, which yields the design matrices $J,K,L$\footnote{Due to their large size, it is not possible to provide these matrices in the paper. Kindly contact the authors for the code and data.}. The result of the state estimation method is demonstrated in Fig.~\ref{fig:state_est}. It can be seen that the observer converges very closely to the true state within three days of the epidemic outbreak. Secondly, Fig.~\ref{fig:pred_y} shows the output prediction/tracking by the observer at every time instant. Notice that, even though it was initialized arbitrarily, the observer converges very fast and minimizes the difference between the measured output and the predicted output.

\textit{Optimal control.}
In the optimal control problem, we choose the cost functional  \eqref{eq:cost_fun} as
\[
J(x,u) = \int_{t_1}^{t_2} \left( x(t)^\T \Gamma x(t) + u(t)^\T \Lambda u(t) \right) \d t
\]
where $\Gamma\geq 0$ and $\Lambda>0$ are chosen to be
\begin{align*}
    \Gamma &= \diag(0.01, 1, 0, 2, 10, 0) \\
    \Lambda &= \diag(0.01, 0.01, 0.01, 0.01).
\end{align*}
In short, we would like to minimize the susceptible, infected, hospitalized, and extinct cases in the time interval $[t_1,t_2]$. The control inputs are bounded as $0\leq u_1(t) \leq 1$, $0\leq u_2(t)\leq 0.9$, $0.1\leq u_3(t)\leq 0.7$, and $0\leq u_4(t)\leq 0.7$. Notice that, in general, it may be physically impossible to make $u_2(t)=1$ because that would mean that inflow to $\E$ is zero, i.e., no one dies due to the epidemic. Moreover, we have $u_3(t)>0$ to allow a baseline diagnosis of the infected cases. If $u_3(t)=0$, then the inflow to $\D$ is zero, i.e., no one is diagnosed with the disease.

The equality constraint in \eqref{eq:constraints} is \eqref{eq:sidher} with estimated parameters (Table~\ref{tab:true_est_parameters}) and initialized at estimated state $\hat{x}(t_1)$. The inequality constraints in \eqref{eq:constraints} are
\[\I(t)-\bar{\I}\leq 0, \; \H(t)-\bar{\H}\leq 0, \; \E(t)-\bar{\E}\leq 0\]
where we choose $\bar{I}=0.5$, $\bar{\H}=0.05$, and $\bar{\E}=0.005$.

\begin{figure}[!thp]
    \centering
    \includegraphics[width=0.475\textwidth]{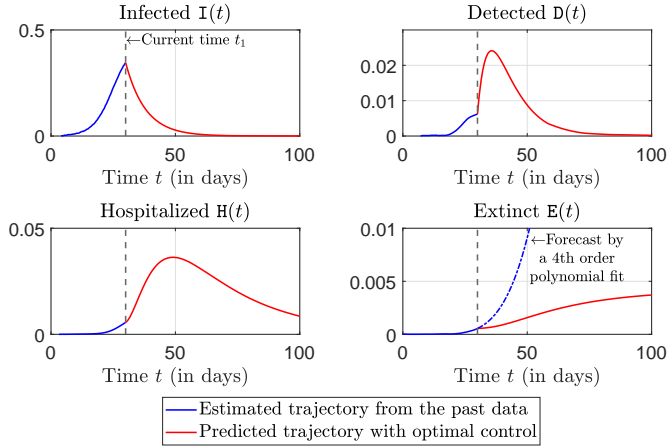}
    \caption{State trajectory with optimal control input.}
    \label{fig:state_optcon}
\end{figure}

\begin{figure}[!thp]
    \centering
    \includegraphics[width=0.475\textwidth]{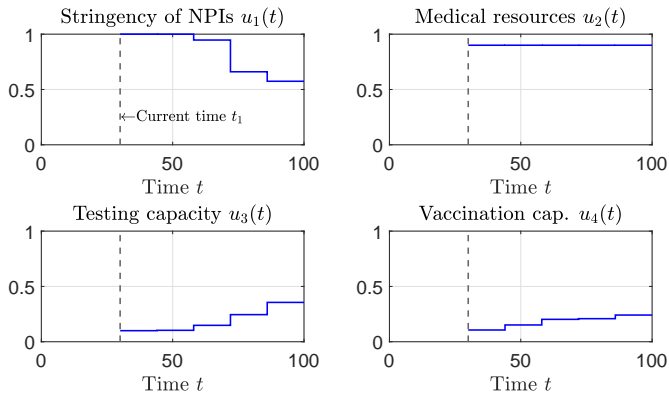}
    \caption{Sequence of optimal control input.}
    \label{fig:optcon_seq}
\end{figure}

We use the \texttt{fmincon} solver in MATLAB for obtaining the optimal solution of \eqref{prob:opt_con}. We consider a piecewise constant control input $u(t)$, where $u(t)$ remains constant for a period of 14 days. The optimal control sequence is obtained for five such periods, i.e., the solution recommends how the policy should be varied in the future after every 14 days. The obtained optimal control input is shown in Fig.~\ref{fig:optcon_seq} and the resulting predicted state trajectory in Fig.~\ref{fig:state_optcon}. Notice from Fig.~\ref{fig:state_optcon} that the constraints on $\I$, $\H$, and $\E$ are satisfied. From the forecast of $\E$, using the polynomial fit of the estimated trajectory, it can be seen that the number of deaths are significantly reduced under the optimal control algorithm. This prediction is reliable only when, in addition to the parameters, the current state $x(t_1)$ is accurately estimated. The results may vary significantly if the unmeasured states are not chosen accurately. From Fig.~\ref{fig:optcon_seq}, we can see that the optimal controller recommends to enact a full lockdown for the first four weeks to suppress the epidemic growth, i.e., to make $\dot{\I}<0$. The lockdown is then lifted gradually by increasing the testing and vaccination capacities. On the other hand, to avoid the number of deaths and to satisfy the hard constraint $\E(t)\leq \bar{\E}$, the optimal controller recommends to employ full medical resources throughout the finite future horizon.

\section{Concluding Remarks}
We presented a unified framework for feedback optimal control of epidemic processes via a general class of nonlinear compartmental models. For the parameter estimation and constrained optimal control, we employed existing methods and techniques from system identification and optimal control theory. For the state estimation, however, the proposed robust observer design criteria is a novel contribution of this paper. By considering a realistic epidemic model, we demonstrated how optimal control policies can be devised by employing the proposed feedback design.

The method proposed in this paper could be helpful in devising optimal control policies during an epidemic outbreak. Robustness and fast convergence of the observer guarantees accurate state estimation of an epidemic model even when few samples of the data are collected. Extending the current method to a moving horizon estimation and control is a straightforward task and will be addressed in the future. Another interesting direction is the observer design for simultaneous parameter and state estimation. Recently proposed learning-based Luenberger observer design method by \cite{niazi2022KKL} could be helpful in this regard. Such an observer is crucial because the existing parameter estimation techniques require an accurate guess of the initial state, which may not be possible in real epidemic scenarios.

\bibliography{bib_ifacconf}

\begin{thebibliography}{40}
\providecommand{\natexlab}[1]{#1}
\providecommand{\url}[1]{\texttt{#1}}
\providecommand{\urlprefix}{URL }
\expandafter\ifx\csname urlstyle\endcsname\relax
  \providecommand{\doi}[1]{doi:\discretionary{}{}{}#1}\else
  \providecommand{\doi}{doi:\discretionary{}{}{}\begingroup
  \urlstyle{rm}\Url}\fi

\bibitem[{Acemoglu et~al.(2021)Acemoglu, Fallah, Giometto, Huttenlocher,
  Ozdaglar, Parise, and Pattathil}]{acemoglu2021}
Acemoglu, D., Fallah, A., Giometto, A., Huttenlocher, D., Ozdaglar, A., Parise,
  F., and Pattathil, S. (2021).
\newblock Optimal adaptive testing for epidemic control: Combining molecular
  and serology tests.
\newblock \emph{arXiv preprint arXiv:2101.00773}.

\bibitem[{Arino and Van~den Driessche(2003)}]{arino2003}
Arino, J. and Van~den Driessche, P. (2003).
\newblock A multi-city epidemic model.
\newblock \emph{Mathematical Population Studies}, 10(3), 175--193.

\bibitem[{Audoly et~al.(2001)Audoly, Bellu, D'Angio, Saccomani, and
  Cobelli}]{audoly2001}
Audoly, S., Bellu, G., D'Angio, L., Saccomani, M.P., and Cobelli, C. (2001).
\newblock Global identifiability of nonlinear models of biological systems.
\newblock \emph{IEEE Transactions on Biomedical Engineering}, 48(1), 55--65.

\bibitem[{Azimi et~al.(2022)Azimi, Sharifi, Fakoorian, Nguyen, and
  Van~Huynh}]{azimi2022}
Azimi, V., Sharifi, M., Fakoorian, S., Nguyen, T., and Van~Huynh, V. (2022).
\newblock State estimation-based robust optimal control of influenza epidemics
  in an interactive human society.
\newblock \emph{Information Sciences}, 592, 340--360.

\bibitem[{Bard(1974)}]{bard1974}
Bard, Y. (1974).
\newblock \emph{Nonlinear Parameter Estimation}.
\newblock Academic Press, New York, USA.

\bibitem[{Barz et~al.(2015)Barz, K{\"o}rkel, Wozny et~al.}]{barz2015}
Barz, T., K{\"o}rkel, S., Wozny, G., et~al. (2015).
\newblock Nonlinear ill-posed problem analysis in model-based parameter
  estimation and experimental design.
\newblock \emph{Computers \& Chemical Engineering}, 77, 24--42.

\bibitem[{Bernard et~al.(2022)Bernard, Andrieu, and Astolfi}]{bernard2022}
Bernard, P., Andrieu, V., and Astolfi, D. (2022).
\newblock Observer design for continuous-time dynamical systems.
\newblock \emph{Annual Reviews in Control}.

\bibitem[{Betts(2010)}]{betts2010}
Betts, J.T. (2010).
\newblock \emph{Practical Methods for Optimal Control and Estimation using
  Nonlinear Programming}.
\newblock SIAM.

\bibitem[{Byrd et~al.(1999)Byrd, Hribar, and Nocedal}]{byrd1999}
Byrd, R.H., Hribar, M.E., and Nocedal, J. (1999).
\newblock An interior point algorithm for large-scale nonlinear programming.
\newblock \emph{SIAM Journal on Optimization}, 9(4), 877--900.

\bibitem[{Byrd et~al.(1987)Byrd, Schnabel, and Shultz}]{byrd1987}
Byrd, R.H., Schnabel, R.B., and Shultz, G.A. (1987).
\newblock A trust region algorithm for nonlinearly constrained optimization.
\newblock \emph{SIAM Journal on Numerical Analysis}, 24(5), 1152--1170.

\bibitem[{Coleman and Li(1994)}]{coleman1994}
Coleman, T.F. and Li, Y. (1994).
\newblock On the convergence of interior-reflective {Newton} methods for
  nonlinear minimization subject to bounds.
\newblock \emph{Mathematical Programming}, 67(1), 189--224.

\bibitem[{Diop and Fliess(1991)}]{diop1991}
Diop, S. and Fliess, M. (1991).
\newblock Nonlinear observability, identifiability, and persistent
  trajectories.
\newblock In \emph{Proc. 30th IEEE Conference on Decision and Control},
  714--719.

\bibitem[{Giordano et~al.(2020)Giordano, Blanchini, Bruno, Colaneri,
  Di~Filippo, Di~Matteo, and Colaneri}]{giordano2020}
Giordano, G., Blanchini, F., Bruno, R., Colaneri, P., Di~Filippo, A.,
  Di~Matteo, A., and Colaneri, M. (2020).
\newblock Modelling the {COVID-19} epidemic and implementation of
  population-wide interventions in {Italy}.
\newblock \emph{Nature Medicine}, 26(6), 855--860.

\bibitem[{Giordano et~al.(2021)Giordano, Colaneri, Di~Filippo, Blanchini,
  Bolzern, De~Nicolao, Sacchi, Colaneri, and Bruno}]{giordano2021}
Giordano, G., Colaneri, M., Di~Filippo, A., Blanchini, F., Bolzern, P.,
  De~Nicolao, G., Sacchi, P., Colaneri, P., and Bruno, R. (2021).
\newblock Modeling vaccination rollouts, {SARS-CoV-2} variants and the
  requirement for non-pharmaceutical interventions in {Italy}.
\newblock \emph{Nature Medicine}, 27(6), 993--998.

\bibitem[{Gomez-Exposito et~al.(2021)Gomez-Exposito, Rosendo-Macias, and
  Gonzalez-Cagigal}]{gomez2021}
Gomez-Exposito, A., Rosendo-Macias, J.A., and Gonzalez-Cagigal, M.A. (2021).
\newblock Monitoring and tracking the evolution of a viral epidemic through
  nonlinear kalman filtering: Application to the {COVID-19} case.
\newblock \emph{IEEE Journal of Biomedical and Health Informatics}.

\bibitem[{Grewal and Glover(1976)}]{grewal1976}
Grewal, M. and Glover, K. (1976).
\newblock Identifiability of linear and nonlinear dynamical systems.
\newblock \emph{IEEE Transactions on Automatic Control}, 21(6), 833--837.

\bibitem[{Hadeler(2011)}]{hadeler2011}
Hadeler, K. (2011).
\newblock Parameter identification in epidemic models.
\newblock \emph{Mathematical Biosciences}, 229(2), 185--189.

\bibitem[{Hermann and Krener(1977)}]{hermann1977}
Hermann, R. and Krener, A. (1977).
\newblock Nonlinear controllability and observability.
\newblock \emph{IEEE Transactions on Automatic Control}, 22(5), 728--740.

\bibitem[{Hethcote(1989)}]{hethcote1989}
Hethcote, H.W. (1989).
\newblock Three basic epidemiological models.
\newblock In \emph{Applied Mathematical Ecology}, 119--144. Springer.

\bibitem[{Hethcote(1994)}]{hethcote1994}
Hethcote, H.W. (1994).
\newblock A thousand and one epidemic models.
\newblock In \emph{Frontiers in Mathematical Biology}, 504--515. Springer.

\bibitem[{Hethcote(2000)}]{hethcote2000}
Hethcote, H.W. (2000).
\newblock The mathematics of infectious diseases.
\newblock \emph{SIAM review}, 42(4), 599--653.

\bibitem[{Hota et~al.(2021)Hota, Godbole, and Par{\'e}}]{hota2021}
Hota, A.R., Godbole, J., and Par{\'e}, P.E. (2021).
\newblock A closed-loop framework for inference, prediction, and control of
  {SIR} epidemics on networks.
\newblock \emph{IEEE Transactions on Network Science and Engineering}, 8(3),
  2262--2278.

\bibitem[{Kirk(2004)}]{kirk2004}
Kirk, D.E. (2004).
\newblock \emph{Optimal Control Theory: An Introduction}.
\newblock Courier Corporation.

\bibitem[{K{\"o}hler et~al.(2021)K{\"o}hler, Schwenkel, Koch, Berberich, Pauli,
  and Allg{\"o}wer}]{kohler2020}
K{\"o}hler, J., Schwenkel, L., Koch, A., Berberich, J., Pauli, P., and
  Allg{\"o}wer, F. (2021).
\newblock Robust and optimal predictive control of the covid-19 outbreak.
\newblock \emph{Annual Reviews in Control}, 51, 525--539.

\bibitem[{Ligon et~al.(2018)Ligon, Fr{\"o}hlich, Chi{\c{s}}, Banga,
  Balsa-Canto, and Hasenauer}]{ligon2018}
Ligon, T.S., Fr{\"o}hlich, F., Chi{\c{s}}, O.T., Banga, J.R., Balsa-Canto, E.,
  and Hasenauer, J. (2018).
\newblock {GenSSI 2.0}: Multi-experiment structural identifiability analysis of
  {SBML} models.
\newblock \emph{Bioinformatics}, 34(8), 1421--1423.

\bibitem[{Ljung(1999)}]{ljung1999}
Ljung, L. (1999).
\newblock \emph{System Identification: Theory for the User}.
\newblock Prentice-Hall, Upper Saddle River, NJ, 2 edition.

\bibitem[{Magal and Webb(2018)}]{magal2018}
Magal, P. and Webb, G. (2018).
\newblock The parameter identification problem for {SIR} epidemic models:
  Identifying unreported cases.
\newblock \emph{Journal of Mathematical Biology}, 77(6), 1629--1648.

\bibitem[{Mei et~al.(2017)Mei, Mohagheghi, Zampieri, and Bullo}]{mei2017}
Mei, W., Mohagheghi, S., Zampieri, S., and Bullo, F. (2017).
\newblock On the dynamics of deterministic epidemic propagation over networks.
\newblock \emph{Annual Reviews in Control}, 44, 116--128.

\bibitem[{Mor{\'e}(1978)}]{more1978}
Mor{\'e}, J.J. (1978).
\newblock The {Levenberg-Marquardt} algorithm: Implementation and theory.
\newblock In \emph{Numerical Analysis}, 105--116. Springer.

\bibitem[{Niazi et~al.(2022)Niazi, Cao, Sun, Das, and Johansson}]{niazi2022KKL}
Niazi, M.U.B., Cao, J., Sun, X., Das, A., and Johansson, K.H. (2022).
\newblock Learning-based design of {Luenberger} observers for autonomous
  nonlinear systems.
\newblock \emph{arXiv preprint arXiv:2210.01476}.

\bibitem[{Niazi and Johansson(2022)}]{niazi2022cdc}
Niazi, M.U.B. and Johansson, K.H. (2022).
\newblock Observer design for the state estimation of epidemic processes.
\newblock \emph{arXiv preprint arXiv:2207.11977}.

\bibitem[{Niazi et~al.(2021)Niazi, Kibangou, Canudas-de Wit, Nikitin, Tumash,
  and Bliman}]{niazi2021arc}
Niazi, M.U.B., Kibangou, A., Canudas-de Wit, C., Nikitin, D., Tumash, L., and
  Bliman, P.A. (2021).
\newblock Modeling and control of epidemics through testing policies.
\newblock \emph{Annual Reviews in Control}, 52, 554--572.

\bibitem[{Nowzari et~al.(2016)Nowzari, Preciado, and Pappas}]{nowzari2016}
Nowzari, C., Preciado, V.M., and Pappas, G.J. (2016).
\newblock Analysis and control of epidemics: A survey of spreading processes on
  complex networks.
\newblock \emph{IEEE Control Systems Magazine}, 36(1), 26--46.

\bibitem[{Par\'{e} et~al.(2020)Par\'{e}, Beck, and Ba\c{s}ar}]{pare2020}
Par\'{e}, P.E., Beck, C.L., and Ba\c{s}ar, T. (2020).
\newblock Modeling, estimation, and analysis of epidemics over networks: An
  overview.
\newblock \emph{Annual Reviews in Control: Special Issue on Systems and Control
  Research Efforts Against COVID-19 and Future Pandemics}, 50, 345--360.

\bibitem[{Par{\'e} et~al.(2018)Par{\'e}, Beck, and Nedi{\'c}}]{pare2018}
Par{\'e}, P.E., Beck, C.L., and Nedi{\'c}, A. (2018).
\newblock Epidemic processes over time-varying networks.
\newblock \emph{IEEE Trans. on Control of Network Systems}, 5(3), 1322--1334.

\bibitem[{Rajaei et~al.(2021)Rajaei, Raeiszadeh, Azimi, and
  Sharifi}]{rajaei2021}
Rajaei, A., Raeiszadeh, M., Azimi, V., and Sharifi, M. (2021).
\newblock State estimation-based control of {COVID-19} epidemic before and
  after vaccine development.
\newblock \emph{Journal of Process Control}, 102, 1--14.

\bibitem[{Saccomani et~al.(2003)Saccomani, Audoly, and
  D'Angi{\`o}}]{saccomani2003}
Saccomani, M.P., Audoly, S., and D'Angi{\`o}, L. (2003).
\newblock Parameter identifiability of nonlinear systems: The role of initial
  conditions.
\newblock \emph{Automatica}, 39(4), 619--632.

\bibitem[{Sethi(2019)}]{sethi2019}
Sethi, S.P. (2019).
\newblock \emph{Optimal Control Theory: Applications to Management Science and
  Economics}.
\newblock Springer Nature Switzerland AG, 3 edition.

\bibitem[{Sontag and Wang(2000)}]{sontag2000}
Sontag, E. and Wang, Y. (2000).
\newblock Lyapunov characterizations of input to output stability.
\newblock \emph{SIAM Journal on Control and Optimization}, 39(1), 226--249.

\bibitem[{Wan et~al.(2007)Wan, Roy, and Saberi}]{wan2007}
Wan, Y., Roy, S., and Saberi, A. (2007).
\newblock Network design problems for controlling virus spread.
\newblock In \emph{Proc. 46th IEEE Conference on Decision and Control},
  3925--3932.

\end{thebibliography}

\end{document}